\documentclass{amsart}
\begin{document}

\vfuzz2pt 
\hfuzz2pt 
\newtheorem{thm}{Theorem}[section]
\newtheorem{corollary}[thm]{Corollary}
\newtheorem{lemma}[thm]{Lemma}
\newtheorem{proposition}[thm]{Proposition}
\newtheorem{defn}[thm]{Definition}
\newtheorem{remark}[thm]{Remark}
\newtheorem{example}[thm]{Example}
\newtheorem{fact}[thm]{Fact}
\
\newcommand{\norm}[1]{\left\Vert#1\right\Vert}
\newcommand{\abs}[1]{\left\vert#1\right\vert}
\newcommand{\set}[1]{\left\{#1\right\}}
\newcommand{\Real}{\mathbb R}
\newcommand{\eps}{\varepsilon}
\newcommand{\To}{\longrightarrow}
\newcommand{\BX}{\mathbf{B}(X)}
\newcommand{\A}{\mathcal{A}}
\newcommand{\onabla}{\overline{\nabla}}
\newcommand{\hnabla}{\hat{\nabla}}
\newcommand{\f}{\mathbf{f}}

\def\proof{\medskip Proof.\ }
\font\lasek=lasy10 \chardef\kwadrat="32 
\def\kwadracik{{\lasek\kwadrat}}
\def\koniec{\hfill\lower 2pt\hbox{\kwadracik}\medskip}

\newcommand*{\C}{\mathbf{C}}
\newcommand*{\R}{\mathbf{R}}
\newcommand*{\Z}{\mathbf {Z}}

\def\sb{f:M\longrightarrow \C ^n}
\def\det{\hbox{\rm det}\, }
\def\detc{\hbox{\rm det }_{\C}}
\def\i{\hbox{\rm i}}
\def\tr{\hbox{\rm tr}\, }
\def\rk{\hbox{\rm rk}\,}
\def\vol{\hbox{\rm vol}\,}
\def\Im {\hbox{\rm Im}\, }
\def\Re{\hbox{\rm Re}\, }
\def\interior{\hbox{\rm int}\, }
\def\e{\hbox{\rm e}}
\def\pu{\partial _u}
\def\pv{\partial _v}
\def\pui{\partial _{u_i}}
\def\puj{\partial _{u_j}}
\def\puk{\partial {u_k}}
\def\div{\hbox{\rm div}\,}
\def\Ric{\hbox{\rm Ric}\,}
\def\r#1{(\ref{#1})}
\def\ker{\hbox{\rm ker}\,}
\def\im{\hbox{\rm im}\, }
\def\I{\hbox{\rm I}\,}
\def\id{\hbox{\rm id}\,}
\def\exp{\hbox{{\rm exp}^{\tilde\nabla}}\.}
\def\cka{{\mathcal C}^{k,a}}
\def\ckplusja{{\mathcal C}^{k+1,a}}
\def\cja{{\mathcal C}^{1,a}}
\def\cda{{\mathcal C}^{2,a}}
\def\cta{{\mathcal C}^{3,a}}
\def\c0a{{\mathcal C}^{0,a}}
\def\f0{{\mathcal F}^{0}}
\def\fnj{{\mathcal F}^{n-1}}
\def\fn{{\mathcal F}^{n}}
\def\fnd{{\mathcal F}^{n-2}}
\def\Hn{{\mathcal H}^n}
\def\Hnj{{\mathcal H}^{n-1}}
\def\emb{\mathcal C^{\infty}_{emb}(M,N)}
\def\M{\mathcal M}
\def\Ef{\mathcal E _f}
\def\Eg{\mathcal E _g}
\def\Nf{\mathcal N _f}
\def\Ng{\mathcal N _g}
\def\Tf{\mathcal T _f}
\def\Tg{\mathcal T _g}
\def\diff{{\mathcal Diff}^{\infty}(M)}
\def\embM{\mathcal C^{\infty}_{emb}(M,M)}
\def\U1f{{\mathcal U}^1 _f}
\def\Uf{{\mathcal U} _f}
\def\Ug{{\mathcal U} _g}
\def\Ric{{\rm {Ric}}}
\def\[f]{{\mathcal U}^1 _{[f]}}
\def\hnu{\hat\nu}
\def\gnu{\nu_g}
\title{Curvature bounded conjugate symmetric statistical structures with complete metric}
\author{Barbara Opozda}




\address{Faculty of Mathematics and Computer Science UJ,
ul. \L ojasiewicza  6, 30-348 Cracow, Poland}

\email{barbara.opozda@im.uj.edu.pl} \maketitle


\begin{abstract}
In the paper two important  theorems about complete affine spheres
are generalized to the case of statistical structures on abstract
manifolds. The assumption about constant sectional curvature is
replaced by the assumption that the curvature satisfies some
inequalities.
\end{abstract}

\section{Introduction}
In the paper we  refer to the following well-known theorems of
affine differential geometry

\begin{thm}{\rm (W. Blaschke, A. Deicke, E. Calabi)}
Let $f:M\to \R^{n+1}$ be an elliptic  affine sphere whose Blaschke
metric is complete. Then  the induced structure on $M$ is trivial.
Consequently the affine sphere is an ellipsoid.
\end{thm}

\begin{thm}{\rm (E. Calabi)}\label{calabi}
 Let $f:M\to\R ^{n+1}$ be a hyperbolic or
parabolic affine sphere whose Blaschke metric is complete. Then the
Ricci tensor of the metric is negative semidefinite.
\end{thm}

The aim of this paper is to generalize these theorems. In
particular, we prove

\begin{thm}\label{inIntroduction}
Let $(g,\nabla) $ be a trace-free conjugate symmetric statistical
structure on a manifold $M$. Assume that $g$ is complete on $M$. If
the sectional $\nabla$-curvature is bounded from below and above on
$M$ then the Ricci tensor of $g$ is bounded from below and above on
$M$. If the sectional $\nabla$-curvature is non-negative everywhere
then the statistical structure is trivial, that is,
$\nabla=\hat\nabla$. If the sectional $\nabla$-curvature is bounded
from $0$ by a positive constant then, additionally,   $M$ is compact
and its first fundamental group is finite.
\end{thm}
 More precise and more general
 formulations of this theorem give Theorems \ref{epsilonfunkcja} and
\ref{thmpositive}. The meaning of the generalization can be
explained as follows. The induced structure on an affine sphere is a
conjugate symmetric trace-free statistical structure. But the
statistical connection on an affine sphere is projectively flat and
its sectional $\nabla$-curvature  is constant. In the theorems we
propose the projective flatness is not needed, which means that the
statistical structure can be non-realizable  on any Blaschke
hypersurface even locally. Moreover, the assumption about the
constant curvature is replaced by the assumption that the curvature
satisfies some inequalities. Since the notion of the sectional
$\nabla$-curvature is relatively new, see  \cite{BW4}, \cite{F}, the
theorems proved in this paper show that the notion is meaningful.

 In the proofs of Theorems
\ref{epsilonfunkcja} and Theorem \ref{thmpositive} we use the same
main tool as in Calabi's theorems, that is, a theorem  on weak
solutions of differential inequalities for the Laplacian of
non-negative functions on complete Riemannian manifolds. In fact, we
shall use only a particular version of this theorem. Note that a
crucial step in the proof of Theorem \ref{epsilonfunkcja} is  an
estimation obtained in Lemma \ref{crucial}. In the case of affine
spheres (Theorem \ref{calabi}) the corresponding part of the proof
is trivial.

\bigskip
\section{Preliminaries}
\subsection{Definitions of statistical structures.}We shall shortly recall basic notions of statistical
geometry. For details we  refer to \cite{BW4}. Let $g$ be a positive
definite Riemannian tensor field on a manifold $M$.
 Denote by $\hnabla$ the
Levi-Civita connection for $g$. A statistical structure is a pair
$(g,\nabla) $ where $\nabla$ is a torsion-free connection such that
the following Codazzi condition is satisfied:

\begin{equation}\label{symmetry}
(\nabla_X g)(Y,Z)=(\nabla _Yg)(X,Z)
\end{equation}
for all $X,Y,Z\in T_x M$, $x\in M$.  A connection $\nabla$
satisfying (\ref{symmetry}) is called a statistical connection for
$g$.

For any connection $\nabla$ one defines its conjugate (dual)
connection $\onabla$ relative to $g$ by the formula
\begin{equation}
g(\nabla _XY,Z)+g(Y,\onabla _XZ)=Xg(Y,Z).
\end{equation}
It is known that the pairs $(g,\nabla)$ and $(g,\onabla)$ are
simultaneously  statistical structures. From now on we assume that
$\nabla$ is a statistical connection for $g$. We have
\begin{equation}\label{R_and_oR}
g(R(X,Y)Z,W)=-g(\overline{R}(X,Y)W,Z),
\end{equation}
where $R$ and $\overline R$ are the curvature tensors for $\nabla$
and $\overline\nabla$, respectively. Denote by $\Ric$ and
$\overline{\Ric}$ the corresponding Ricci tensors. Note that, in
general, these Ricci tensors are not necessarily symmetric. The
curvature and the Ricci tensor of $\hnabla$ will be denoted by $\hat
R$ and $\widehat{\Ric}$, respectively.
 The
function
\begin{equation}
\rho=\tr _g\Ric(\cdot,\cdot)
\end{equation}
is the scalar curvature of $(g,\nabla)$. Similarly, one can define
the scalar curvature $\overline\rho$ for $(g,\onabla)$ but, by
(\ref{R_and_oR}), $\rho=\overline\rho$. We also have the usual
scalar curvature $\hat\rho$ for $g$.

Denote by $K$ the difference tensor between $\nabla$ and $\hnabla$,
that is,
\begin{equation}
\nabla _XY=\hnabla _XY+K_XY.
\end{equation}
Then
\begin{equation}
\onabla_XY=\hnabla _XY-K_XY.
\end{equation}
 $K(X,Y)$ will stand for $K_XY$. Since $\nabla$  and $\hnabla$ are  without
torsion, $K$  as a $(1,2)$-tensor is symmetric. We have $ (\nabla
_Xg)(Y,Z)=(K _Xg)(Y,Z)=-g(K_XY,Z)-g(Y,K_XZ)$. It is now clear that
the symmetry of $\nabla g$ and $K$ implies the symmetry of $K_X$
relative to $g$ for each $X$. The converse also holds. Namely, if
$K_X$ is symmetric relative to $g$ then we have $( \nabla _Xg)(Y,Z)
= -2g(K_XY,Z)$.

We define the statistical cubic  form $A$ by
\begin{equation}
A(X,Y,Z)=g(K_XY,Z).
\end{equation}
It is clear that a statistical structure can be defined equivalently
as a pair $(g,K)$, where $K$ is a symmetric tensor field of type
$(1,2)$ which is also symmetric relative to $g$, or  as a pair
$(g,A)$, where $A$ is a symmetric cubic form.

A statistical structure is trace-free if $\tr_gK(\cdot,\cdot)=0$
(equivalently, $\tr_gA(X,\cdot,\cdot)=0$ for every $X$;
equivalently, $\tr K_X=0$ for every $X$). The trace-freeness is also
equivalent to the condition that $\nabla \nu_g=0$, where $\nu_g$ is
the volume form determined by $g$. In affine differential geometry
the trace-freeness is called the apolarity. The assumption  about
the trace-freeness  of a statistical structure is essential in  all
the theorems mentioned in the Introduction.
\medskip

\subsection{Relations between curvature tensors of  statistical structures.}
It is known that
\begin{equation}\label{from_Nomizu_Sasaki}
R(X,Y)=\hat R(X,Y) +(\hnabla_XK)_Y-(\hnabla_YK)_X+[K_X,K_Y].
\end{equation}
Writing the same equality for $\onabla$ and adding both equalities
we get
\begin{equation}\label{R+oR}
R(X,Y)+\overline R(X,Y) =2\hat R(X,Y) +2[K_X,K_Y].
\end{equation}
In particular, if $R=\overline R$ then
\begin{equation}
R(X,Y)=\hat R(X,Y) +[K_X,K_Y],
\end{equation}
which  can be shortly written as
\begin{equation}\label{short}
R=\hat R +[K,K].
\end{equation}
 Using (\ref{R+oR}) and assuming that  a given
statistical structure is trace-free one gets, see \cite{BW4},
\begin{equation}\label{Ric+oRic}
\Ric (Y,Z)+\overline{\Ric} (Y,Z)=2\widehat{\Ric}(Y,Z) -2g(K_Y,K_Z).
\end{equation}
In particular, if $(g, \nabla )$ is trace-free then
\begin{equation}\label{Riccihat_ge}
2\widehat{\Ric}(X,X)\ge
\Ric(X,X)+\overline{\Ric}(X,X).\end{equation} If, moreover,
$R=\overline R$ then
\begin{equation}\label{for R=oR}
\widehat{\Ric}\ge \Ric.
\end{equation}
 The following lemma
follows from  formulas (\ref{R_and_oR})
and(\ref{from_Nomizu_Sasaki}).

\begin{lemma}\label{przeniesiony_lemat}
Let $(g,\nabla) $ be a statistical structure. The following
conditions are equivalent:
\newline
{\rm 1)} $R=\overline R$,
\newline
{\rm 2)} $\hnabla K$ is symmetric (equiv. $\hat\nabla A$ is
symmetric),
\newline
{\rm 3)} $g(R(X,Y)Z,W)$ is skew-symmetric relative to $Z,W$.
\end{lemma}
The family of statistical  structures satisfying one of the above
conditions is as important in the theory of statistical structures
as the family of affine spheres in affine differential geometry.  A
statistical structure satisfying 2) in the above lemma was called in
\cite{L} conjugate symmetric. We shall adopt this definition.

Note that the condition $R=\overline R$ implies the symmetry of
$\Ric$.

Taking now the trace  relative to $g$ on both sides of
(\ref{Ric+oRic}) and taking into account that $\rho =\overline\rho$,
we get
\begin{equation}\label{theorema_egregium}
\hat\rho =\rho +\Vert K\Vert^2=\rho +\Vert A\Vert ^2
\end{equation}
for a trace-free statistical structure.
\medskip

\subsection{Sectional $\nabla$-curvature.}
The notion of a sectional $\nabla$-curvature was first introduced in
\cite{BW4}.  Namely, the tensor field
\begin{equation}
\mathcal R= \frac{1}{2}(R+\overline R)
\end{equation}
is a Riemannian curvature tensor. In particular,  it satisfies the
condition $$g(\mathcal R(X,Y)Z,W)=-g(\mathcal R(X,Y)W,Z).$$  In
general, this condition is not satisfied by the curvature tensor
$R$. In the case where a given statistical structure is conjugate
symmetric the curvature tensor $R$ satisfies this condition. In
\cite{BW4} we defined the sectional $\nabla$-curvature by
\begin{equation}
k(\pi) =g(\mathcal R(e_1,e_2)e_2,e_1)
\end{equation}
for a vector plane $\pi \in T_xM$, $x\in M$ and $e_1,e_2$ any
orthonormal basis of $\pi$. It is a well-defined notion.

In general, Schur's lemma does not hold for the sectional
$\nabla$-curvature. But, if a statistical structure is conjugate
symmetric (in this case $\mathcal R=R$) then the second Bianchi
identity holds and, consequently, the Schur lemma holds, see
\cite{BW4}. Thus, in Theorems 3.4 and 4.1 from \cite{No} (which were
the inspiration for our investigation), the functions $\tilde R$
are, in fact, constant if $n>2$.

\medskip

\subsection{Statistical structures on  affine hypersurfaces.}
 The theory of affine hypersurfaces in $\R ^{n+1}$ is a natural
source of statistical structures. For the theory we refer to
\cite{LSZ} or \cite{NS}. We recall here only some selected facts.

 Let $\mathbf{f} :M\to \R
^{n+1}$ be a locally strongly convex  hypersurface.  For simplicity
assume that $M$ is connected and orientable. Let $\xi$ be a
transversal vector field on $M$. We have the induced  volume form
$\nu _\xi$ on $M$  defined as follows
$$\nu_\xi (X_1,...,X_n)=\det (\mathbf {f}_*X_1,...,\mathbf{f}_*X_n,\xi).$$
We  also have the induced connection $\nabla$ and the second
fundamental form $g$ defined by the Gauss formula
$$D_X\mathbf{f}_*Y=\mathbf{f}_*\nabla _XY +g(X,Y)\xi,$$
where $D$ is the standard flat connection on $\R ^{n+1}$. Since the
hypersurface is locally strongly convex, the second fundamental form
$g$ is definite. By multiplying $\xi$ by $-1$ if necessary, we can
assume that $g$ is positive definite. A transversal vector field is
called equiaffine if $\nabla \nu_\xi=0$. This condition is
equivalent to the fact that $\nabla g$ is symmetric, i.e.
$(g,\nabla)$ is a statistical structure. It means, in particular,
that for a statistical structure obtained on a hypersurface by a
choice of an equiaffine  transversal vector field, the Ricci tensor
of $\nabla$ is automatically symmetric. A hypersurface equipped with
an equiaffine transversal vector field and the induced structure is
called an equiaffine hypersurface.

Recall now the notion of the shape operator. Having a fixed
equiaffine transversal vector field $\xi$  and differentiating it we
get the Weingarten formula
$$D_X\xi= -\mathbf{f}_*\mathcal SX.$$
The tensor field $\mathcal S$ is called the shape operator for
$\xi$. If $R$ is the curvature tensor for the induced connection
$\nabla$ then
\begin{equation}\label{Gauss_equation_for R}
R(X,Y)Z=g(Y,Z)\mathcal SX-g(X,Z)\mathcal SY.
\end{equation}
This is the Gauss equation for $R$. The Gauss equation for the dual
structure is the following
\begin{equation}\label{Gauss_equation_for_oR}
\overline R(X,Y)Z=g(Y,\mathcal SZ)X-g(X,\mathcal SZ)Y.
\end{equation}
It follows that the dual connection is projectively flat if $n>2$.
The dual connection is also projectively flat for 2-dimensional
surfaces equipped with an equiaffine transversal vector field, that
is,    $\overline \nabla \overline{\Ric}$ is symmetric. The form
$g(\mathcal SX,Y)$ is symmetric for any equiaffine transversal
vector field.

We have the volume form $\nu_g$ determined by $g$ on $M$. In
general, this volume form is not covariant constant relative to
$\nabla$. The starting point of the classical affine differential
geometry is the theorem saying that there is a unique equiaffine
transversal vector field $\xi$ such that $\nu_\xi =\nu_g$. This
unique transversal vector field is called the affine normal vector
field or the Blaschke affine normal. The second fundamental form for
the affine normal is called the Blaschke metric. A non-degenerate
hypersurface endowed with the affine Blaschke normal is called a
Blaschke hypersurface. The induced statistical structure is
trace-free on a Blaschke hypersurface. If the affine lines
determined by the affine normal vector field meet at one point or
are parallel then the hypersurface is called an affine sphere. In
the first case the sphere is called proper in the second one
improper. The class of affine spheres is very large. There exist
many conditions characterizing affine spheres. For instance, a
Blaschke hypersurface is an affine sphere if and only if
$R=\overline R$. Therefore, conjugate symmetric statistical
manifolds  can be regarded as generalizations of affine spheres. For
connected affine spheres the shape operator $S$ is a constant
multiple of the identity, i.e., $S=\kappa\, \id$ for some constant
$\kappa$.

 If we choose a positive definite Blaschke
metric on a connected locally strongly convex affine sphere then we
call the sphere elliptic if $\kappa >0$, parabolic if $\kappa =0$
and hyperbolic if $\kappa <0$.
\medskip

\subsection{Conjugate symmetric statistical structures
non-realizable on affine spheres} As we have already mentioned, if
$\nabla$ is a connection on a hypersurface induced by an equiaffine
transversal vector field then the conjugate connection $\onabla$ is
projectively flat. Therefore the projective flatness of the
conjugate connection is a necessary condition for $(g,\nabla)$ to be
realizable as the induced structure on a hypersurface equipped with
an equiaffine transversal vector field. In fact, one of the
fundamental theorems in affine differential geometry (see, e.g.
\cite{NS}) says, roughly speaking,  that it is also a sufficient
condition for the local realizability of a Ricci-symmetric
statistical structure, but we will not need it in this paper. Note
also that, if $(g,\nabla)$ is a conjugate symmetric statistical
structure then $\nabla$ and $\overline \nabla$ are simultaneously
projectively flat. Indeed, it is obvious for $n>2$. If $n=2$ we can
argue as follows. It suffices to prove that if $\overline\nabla$ is
projectively flat then so is $\nabla$. Since $R=\overline R$,
$\nabla$ is Ricci-symmetric. By the fundamental theorem mentioned
above $(g,\nabla)$ can be locally realized on an equiaffine surface
in $\R ^3$. By Lemma 12.5 from \cite{BW4} the surface is an
equiaffine sphere, that is, the shape operator is locally a constant
multiple of the identity, and hence $\nabla$ is projectively flat.
It follows that if $(g,\nabla)$ is conjugate symmetric then it is
locally realizable on an equiaffine hypersurface if only if $\nabla$
or $\overline\nabla$ is projectively flat.


We shall now consider trace-free conjugate symmetric statistical
structures.
 The following fact was observed in \cite{BW4}, see Proposition 4.1 there. If
$(g,\nabla)$ is the induced statistical structure on an affine
sphere, the metric $g$ is not of constant sectional curvature and
$\alpha\ne 1, -1$ is a real number,  then $\nabla^\alpha:=\hnabla
+\alpha K$ is not projectively flat and therefore it cannot be
realized (even locally) on any affine sphere. Of course,
$(g,\nabla^\alpha)$ is again a statistical conjugate symmetric
structure (by 2) of Lemma \ref{przeniesiony_lemat}) and since the
initial structure was trace-free (because an affine  sphere is
endowed with the Blaschke structure),   $(g,\nabla^\alpha)$ is
trace-free as well. Note also that there are very few affine spheres
whose Blaschke metric has constant sectional curvature, see
\cite{LSZ}, which means that the assumption that $g$ is not of
constant sectional curvature is not restrictive.

The following example shows another easy way of producing conjugate
symmetric trace-free statistical structures which are non-realizable
(even locally) on affine spheres.

Let $M=\R^n$, where $n\ge 4$, be equipped with the standard  flat
metric tensor field $g$. Let $x^1,...,x^n$ be the canonical
coordinate system and $e_1,...,e_n$ be the canonical orthonormal
frame. Define the cubic form $A=(A_{ijk})$ on $M$, where
$A_{ijk}=A(e_i,e_j,e_k)$, by
\begin{equation}\label{w_przykladzie}
\begin{array}{rcl}
&&A_{ijk}=0 \  \ if \ at \ least\ two\ of\ indices \ i,j,k \ are\
equal,\\
&& A_{ijk}\in \R^+\ if\ the\ indices\ i,j,k\ are\ mutually \
distinct.
\end{array}
\end{equation}
Then $\hat\nabla K=0$ and, consequently, $R=\overline R$. Observe
now that the  connection $\overline\nabla=\hnabla -K$, where
$g(K(X,Y),Z)=A(X,Y,Z)$, is not projectively flat and therefore,
$(g,\nabla)$ cannot be  realized on any Blaschke hypersurface, even
locally. Indeed, suppose that $\overline\nabla$ is projectively
flat. Then we must have $g(R(e_i,e_j)e_j, e_l)=0$ for $i\ne j$ and
$l\ne i,j$. On the other hand, by (\ref{from_Nomizu_Sasaki}), we
have
\begin{eqnarray*}
&&g(R(e_i,e_j)e_j,e_l)=g([K_{e_i},K_{e_j}]e_j,e_l)=-g(K_{e_j}K_{e_i}e_j,e_l)\\
&&\ \ \ \ =-g(K(e_i,e_j),K(e_j,e_l))=-\sum_{s=1}^nA_{ijs}A_{jls}.
\end{eqnarray*}
By  (\ref{w_przykladzie}) it is clear that the function
$-\sum_{s=1}^nA_{ijs}A_{jls}$ is negative if $n\ge 4$.

Another version of this example (with $\hat\nabla K\ne 0$) is given
by the symmetric $A_{ijk}$, where
\begin{equation}
\begin{array}{rcl}
&&A_{ijk}=0 \  \ if \ at \ least\ two\ of\ indices \ i,j,k \ are\
equal,\\
&& (A_{ijk})_x=
x^1+...+\widehat{(x^i)}+...+\widehat{(x^j)}+...+\widehat{(x^k)}+...+x^n
\  for\  i<j<k,
\end{array}
\end{equation}
where $\widehat{(x^l)}$ means that the coordinate $x^l$ was removed
from the sum. One  can easily check that $\hat\nabla A$ is
symmetric. Indeed, we want to check that $\partial_lA_{ijk}=\partial
_iA_{ljk}$ for $l\ne i$. It is sufficient to assume that $j\ne k$.
Consider the cases: a) $l=j$ or $l=k$, b) $i=j$ or $i=k$, $\l\ne j
$, $l\ne k$, c) $i\ne j$, $i\ne k$, $l\ne j$, $l\ne k$. In cases a)
and b) both sides of the required equality vanish. In the last case,
where all indices are mutually distinct, on both sides of the
required equality we get 1.

In the same manner as in the previous example  one  sees  that
$\nabla$ is not projectively flat on $(\R^+)^n$ if $n\ge 4$.

The  considerations of this subsection show that the class of
conjugate symmetric trace-free statistical  manifolds is much larger
than the class of affine spheres, even in the local setting.


\bigskip

\section{Curvature bounded conjugate symmetric trace-free statistical structures }
Let $n=\dim M$ and $(g,\nabla)$ be a statistical structure on $M$.
From now on we assume that the structure is trace-free and conjugate
symmetric.
 Assume moreover that
  \begin{equation}\label{warunek12}
  H_2\le k(\pi)\le H_1,\end{equation}
  for every vector plane $\pi\subset T_xM$ and $x\in M$.
   Denote by $\varepsilon$ the
difference $H_1-H_2$ and set
$$H_3=H_2-\frac{n-2}{2}\varepsilon.$$
The quantities $H_1$, $H_2$ and $\varepsilon$ can be functions on
$M$ (not satisfying any smoothness assumptions), but in the main
theorem of this section, that is, in Theorem \ref{epsilonfunkcja},
$H_3$ must be a real number. The condition (\ref{warunek12}) can be
written as
\begin{equation}\label{warunek1}
H_1-\varepsilon\le k(\pi)\le H_1
\end{equation}
or
\begin{equation}\label{warunek3}
H_3+\frac {n-2}{2}\varepsilon\le k(\pi)\le
H_3+\frac{n}{2}\varepsilon.
\end{equation}

\begin{thm}\label{epsilonfunkcja}
Let $(g,\nabla)$ be a trace-free conjugate symmetric statistical
structure on an $n$-dimensional  manifold $M$. Assume that $(M,g)$
is complete and the sectional $\nabla$-curvature   $k$ satisfies the
inequalities {\rm(\ref{warunek3})} on $M$,  where $H_3$ is a
non-positive number and $\varepsilon$ is a non-negative function on
$M$. Then the Ricci tensor $\widehat \Ric$ of $g$ satisfies the
inequalities
\begin{equation}\label{main_inequality}
(n-1)H_3+\frac{(n-1)(n-2)}{2}\varepsilon\le\widehat\Ric\le-(n-1)^2H_3
+\frac{(n-1)n}{2}\varepsilon.
\end{equation}
The scalar curvature $\hat\rho$ of $g$ satisfies the inequalities
\begin{equation}
n(n-1)H_3+\frac{(n-1)(n-2)n}{2}\varepsilon\le\hat\rho\le\frac{n^2(n-1)}{2}\varepsilon.
\end{equation}
\end{thm}
 \proof In what follows the scalar multiplication $g$  will be also
denoted by $\langle\ , \ \rangle$.  The following lemma is crucial
in the following proof
\begin{lemma}\label{crucial}
Let $V$ be any unit vector  of $T_pM$.  Denote by $T_V$ the
$(0,4)$-tensor given by
\begin{equation}
T_V(X,Y,Z,W)=-\langle K_VX,R(Y,Z)W\rangle-2\langle
K_VW,R(Y,Z)X\rangle.
\end{equation}Assume that
\begin{equation}
H_3+\frac{n-2}{2}\varepsilon\le k( \pi)\le
H_3+\frac{n}{2}\varepsilon
\end{equation}
 for
some $H_3\in \R$, $\varepsilon\in \R^+$ and for all vector planes
$\pi\subset T_pM$. Then
\begin{equation}
\langle T'_V,A_V\rangle\ge (n+1)H_3 \psi_V,
\end{equation}
where
\begin{equation}
A_V(X,Z)=A(V,X,Z),
\end{equation}
\begin{equation}
T'_V(X,Z)=\tr _gT_V(X,\cdot,Z,\cdot)
\end{equation}
and
\begin{equation}
\psi _V=\langle A_V,A_V\rangle.
\end{equation}
\end{lemma}

Proof of  Lemma \ref{crucial}. Let $e_1,...,e_n$ be an eigenbasis of
$K_V$ and $K_V e_i=\lambda_ie_i$ for $i=1,..., n$. Then
$\psi_V=\lambda_1^2+...+\lambda_n^2$. We have

\begin{eqnarray*}
&&\langle {T'} _V, A_V\rangle = -\sum_{i,j,k}[\langle
K_Ve_j,R(e_i,e_k)e_i\rangle\langle K_Ve_j,e_k\rangle+2\langle
K_Ve_i, R(e_i,e_k)e_j\rangle\langle
K_Ve_j,e_k\rangle]\\
&&\ \ \ \ \ \ \ \ \ \ \ \ \ \ =
\sum_{i,j}(\lambda_j^2-2\lambda_i\lambda_j) k_{ij},
\end{eqnarray*}
where $k_{ij}=k(e_i\wedge e_j)$. Since $k_{ij}=k_{ji}$ and
$k_{ii}=0$, we obtain
\begin{equation}\label{najpodstawowsza}
\begin{array}{rcl}
&&\langle T'_V,A_V \rangle
=(\lambda_1^2k_{11}+...+\lambda_1^2k_{1n})+...+(\lambda_n^2k_{n1}+...+\lambda_n^2k_{nn})
-4\sum_{i<j}\lambda_i\lambda_jk_{ij}\\
&&\ \ \ \ \
 \ \ \ \ \ \ \ \ =\sum_{i<j}(\lambda_j-\lambda_i)^2k_{ij}-2\sum_{i<j}\lambda_i\lambda_j
k_{ij}.
\end{array}
\end{equation}
In the last term we now replace $\lambda_n$ by $-\lambda
_1-...-\lambda _{n-1}$. We get

\begin{eqnarray*}
&&-\sum_{i<j}\lambda _i\lambda_jk_{ij}=-\lambda_1\lambda_2k_{12}-\ \
\ \ ...\ \ \
-\lambda_1(-\lambda_1-...-\lambda_{n-1})k_{1n}\\
&&\ \ \ \ \ \  \ \ \  \ \ \ \ \ \ \ \ \ \ \ \ \ \ \ \ \ \ \
-\lambda_2\lambda_3k_{23}-...-\lambda_2(-\lambda
_1-...-\lambda _{n-1})k_{2n}\\
&&\ \ \ \ \  \ \ \ \ \ \ \ \ \ \ \ \ \ ...\\
&&\ \ \ \ \ \ \ \  \ \ \ \ \ \ \ \ \ \ \ \  \ \ \ \ \ \ \ \ \ \ \ \
\ \ \ \ \ \ \ \
-\lambda_{n-1}(-\lambda _1-...-\lambda _{n-1})k_{n-1,n}\\
&&\ \ \ \ \ \ \ \ \ \ \ \  \ \ \ \ \  \ \ =
 \lambda _1\lambda_2(-k_{12}+k_{1n})+...+\lambda_1\lambda_{n-1}(-k_{1,n-1}+k_{1n})+\lambda_1^2k_{1n}\\
&&\ \ \ \ \ \ \  \
+\lambda_1\lambda_2k_{2n}+\lambda_2\lambda_3(-k_{23}+k_{2n})+...+\lambda_2\lambda_{n-1}
(-k_{2,n-1}+k_{2n})+\lambda_2^2k_{2n}\\
&&\ \ \ \ \ \ \ \ \ \ \ \ \ \ \ \ \ \ \  \ \ \ \ \ \ \  \ \ \ \ \ \
... \\
&& \ \ \ \ \ \ \  \ \ \ \ \ \ \ \ \ \ \ \ \ \ \ \ \ \ \ \
 +\lambda_{n-1}\lambda_1k_{n-1,n}+...+\lambda_{n-1}\lambda_{n-2}k_{n-1,n}+
\lambda_{n-1}^2k_{n-1,n}\\
&&\ \ \ \ \ \ \ \ \ \ \ \ \  \ \ \ \ \ \ \ =\sum_{i<j\le
n-1}\lambda_i\lambda_j(k_{in}+k_{jn}-k_{ij})+\sum_{i=1}^{n-1}\lambda
_i^2k_{in}.
\end{eqnarray*}
Thus, using the  assumption (\ref{warunek12}) and the condition
$\lambda_n=-\lambda_1-...-\lambda_{n-1}$, we get

\begin{eqnarray*}
\langle T'_VA_V\rangle  &&\ge \sum_{i<j\le
n}(\lambda_i-\lambda_j)^2H_2
+2\sum_{i=1}^{n-1}\lambda_i^2H_2+2\sum_{i<j\le
n-1}\lambda_i\lambda_j(k_{in}+k_{jn}-k_{ij})\\
&&=\sum_{i<j\le
n-1}(\lambda_i^2+\lambda_j^2-2\lambda_i\lambda_j)H_2+\sum_{i=1}^{n-1}(\lambda_i-\lambda_n)^2H_2\\
&&\ \ \ \ \ \ \ \ \ \ \ \
 \ \ \ \  \ \ \ \ \ \ \ \ \ \ \ \ \ \ \ \ \ \ +2\sum_{i=1}^{n-1}\lambda_i^2H_2+2\sum_{i<j\le
n-1}\lambda_i\lambda_j(k_{in}+k_{jn}-k_{ij})\\
&&=\sum_{j<j\le n-1}(\lambda_i^2+\lambda_j^2)H_2 +2\sum _{i<j\le
n-1}\lambda _i\lambda_j (k_{in} +k_{jn}-k_{ij}-H_2)
\\
&&\ \ \ \ \ \ \ \ \ \ \ \ \ \ \ +\sum_{i=1}^{n-1}\lambda_i^2H_2
+(n-1)\lambda_n^2H_2-2\sum_{i=1}^{n-1} \lambda_i\lambda_nH_2
+2\sum_{i=1}^{n-1}\lambda_i^2H_2\\
&&= (n+1)\sum_{i=1}^{n-1}\lambda_i^2H_2 +2\sum _{i<j\le n-1}\lambda
_i\lambda_j (k_{in} +k_{jn}-k_{ij}-H_2)+(n+1)\lambda_n^2H_2\\
&& =(n+1)\psi_VH_2 +2\sum _{i<j\le n-1}\lambda _i\lambda_j (k_{in}
+k_{jn}-k_{ij}-H_2).
\end{eqnarray*}
Therefore, it is sufficient to prove

\begin{equation}\label{podstawowa}
(n+1)\psi_V(H_2-H_3) +2\sum _{i<j\le n-1}\lambda _i\lambda_j (k_{in}
+k_{jn}-k_{ij}-H_2)\ge 0.
\end{equation}
The left hand side of this inequality can be  written and then
estimated as follows
\begin{equation}\label{ostateczne}
\begin{array}{rcl}
&&(n+1)(\lambda_1^2+...+\lambda
_{n-1}^2)(H_2-H_3)+n\lambda_n^2(H_2-H_3)\\&&\ \ \ \ \ \ \ \ \ \ \ \
\ \ \ \
\ \  \ \ \ \ \ \ \  +(\lambda_1+...+\lambda_{n-1})^2(H_2-H_3)\\
&&\ \ \ \ \ \ \ \ \ \ \ \ \ \ \ \ \ \ \ \ \ \  \ +2\sum
_{i<j\le n-1}\lambda _i\lambda_j (k_{in} +k_{jn}-k_{ij}-H_2)\\
&&\ge\frac{n+2}{n-2}(H_2-H_3)(n-2)(\lambda_1^2+...+\lambda_{n-1}^2)\\
&&\ \ \ \ \ \ \ \ \ \ \ \ \ \ \ \ \ \ \ \ \ \ \ \ \ \ \
+2\sum_{i<j\le n-1}\lambda_i\lambda_j(k_{in}+k_{jn}-k_{ij}-H_3).
\end{array}
\end{equation}
Assume now that $n\ge 4$. Observe that for $i<j\le n-1$ we have

$$k_{in}+k_{jn}-k_{ij}-H_3\ge 0.$$
Indeed, we have $ k_{in}+k_{jn}-k_{ij}-H_3\ge
2H_2-H_1-H_3=(\frac{n}{2}-2)\varepsilon\ge 0$ for $n\ge 4$. Moreover
$$ \frac{n+2}{n-2}(H_2-H_3)\ge k_{in}+k_{jn}-k_{ij}-H_3.$$
 Namely, since $H_1=H_3+\frac{n}{2}\varepsilon$ and
$H_2=H_3+\frac{n-2}{2}\varepsilon$, we have

\begin{eqnarray*}
k_{in}+k_{jn}-k_{ij}-H_3\le
2H_1-H_2-H_3=\frac{n+2}{2}\varepsilon=\left(\frac{n+2}{n-2}\right)
\left(\frac{n-2}{2}\varepsilon\right)=\frac{n+2}{n-2}(H_2-H_3).
\end{eqnarray*}
We now can make father estimations in (\ref{ostateczne}) as follows

\begin{eqnarray*}
&&\frac{n+2}{n-2}(H_2-H_3)(n-2)(\lambda_1^2+...+\lambda_{n-1}^2)+2\sum_{i<j\le
n-1}\lambda_i\lambda_j(k_{in}+k_{jn}-k_{ij}-H_3)\\
&&\ge
(n-2)(\lambda_1^2+...+\lambda_{n-1}^2)(k_{in}+k_{jn}-k_{ij}-H_3)+2\sum
_{i<j\le
n-1}\lambda_i\lambda_j(k_{in}+k_{jn}-k_{ij}-H_3)\\
&&=\sum_{i<j\le
n-1}(\lambda_i+\lambda_j)^2(k_{in}+k_{jn}-k_{ij}-H_3)\ge 0.
\end{eqnarray*}
The lemma is proved for $n\ge 4$. Consider now the case $n=3$. By
the trace-freeness we can assume that  $\lambda_1\lambda_2\ge 0$. We
compute and estimate the left hand side of (\ref{podstawowa}) as
follows
\begin{eqnarray*}
&&2(\lambda _1^2+\lambda_2^2+\lambda_3^2)\varepsilon
+2\lambda_1\lambda_2(k_{13}+k_{23}-k_{12}-H_2)\\
&&\ge 2(\lambda _1^2+\lambda_2^2+\lambda_3^2)\varepsilon
+2\lambda_1\lambda_2(2H_2-H_1-H_2)\\
&&=2(\lambda _1^2+\lambda_2^2+\lambda_3^2)\varepsilon
-2\lambda_1\lambda_2\varepsilon\\
&& =(\lambda_1-\lambda_2)^2\varepsilon +(\lambda_1^2
+\lambda_2^2+\lambda_3^2)\varepsilon\ge 0.
\end{eqnarray*}
Finally, consider the case $n=2$. In this case we have $H_2=H_3$,
$\lambda_2=-\lambda_1$ and $\psi_V=2\lambda_1^2$. Going back to
(\ref{najpodstawowsza}) we get
$$\langle T'_V,A_V\rangle=6\lambda_1^2k_{12}\ge 3\psi_VH_3.$$
The proof of  Lemma \ref{crucial} is completed. \koniec

It is well-known that for any tensor field $s$ the following formula
holds
\begin{equation}\label{Simons}
\Delta (g(s,s))=2g(\Delta s,s) +2g(\hat\nabla s,\hat\nabla s),
\end{equation}
where $\Delta s$ is defined by
\begin{equation}
\Delta s =\sum _{i=1}^n \hat\nabla^2_{e_ie_i}s,
\end{equation}
for any orthonormal frame $e_i$. We shall now compute $\Delta \psi$
for
\begin{equation}
\psi=g(A,A).
\end{equation}
Let $p\in M$, $X,Y,Z\in T_pM$ and $e_1,...,e_n$ be an orthonormal
basis of $T_pM$. Extend all these vectors  along
$\hat\nabla$-geodesics starting at $p$ and denote the obtained
vector fields by the same letters $X,Y,X$, $e_1,...,e_n$,
respectively. Of course, $\hat\nabla X=\hat\nabla Y=\hat\nabla
Z=0$,$\hat\nabla e_1=0$,..., $\hat\nabla e_n=0$  at $p$. The frame
field $e_1,..., e_n$ is orthonormal. Since $\hnabla A$ is symmetric,
one gets at $p$

\begin{eqnarray*}
&&\sum_{i=1}^n(\hat\nabla_{e_ie_i}^2A)(X,Y,Z)=
\sum_{i=1}^n(\hnabla_{e_i}(\hnabla A))(e_i,X, Y,Z)=
\sum_{i=1}^n\hnabla_{e_i}((\hnabla _{e_i}A)(X,Y,Z))\\
&& \ \ \ \ \ \ \ =\sum_{i=1}^n\hnabla_{e_i}
((\hnabla _{X}A)(e_i,Y,Z))=\sum_{i=1}^n(\hnabla_{e_i}(\hnabla _{X}A))(e_i,Y,Z)) \\
&&\ \ \ \ \ \ \ =\sum_{i=1}^n(\hat R (e_i,X) A)(e_i,Y,Z)+
\sum_{i=1}^n(\hnabla_{X}(\hnabla _{e_i}A))(e_i,Y,Z))\\
&&\ \ \ \ \ \ \   =\sum_{i=1}^n(\hat R (e_i,X) A)(e_i,Y,Z)+
\sum_{i=1}^n\hnabla_{X}((\hnabla _{e_i}A)(e_i,Y,Z)).
\end{eqnarray*}
Thus
\begin{equation}
(\Delta A)(X,Y,Z)= \tr _g (\hat R(\cdot, X)A)(\cdot, Y,Z).
\end{equation}
Since $\hat R =R -[K,K]$,
 we have
\begin{equation}\label{deltaA}
(\Delta A)(X,Y,Z)= \tr_g (R(\cdot,X)A)(\cdot,
Y,Z)-\tr_g([K_{\cdot},K_X]A)(\cdot. Y,Z).
\end{equation}
We shall use the following inequality proved, in fact,  on p. 84 in
\cite{NS}.
\begin{proposition}\label{nomizu}
For a trace-free statistical structure we have
\begin{equation}
g(F,A)\ge\frac{n+1}{n(n-1)}(g(A,A))^2,
\end{equation}
where
\begin{equation}
F(X,Y,Z)= -\tr _ g([K_{\cdot},K_X]A)(\cdot , Y,Z).
\end{equation}
\end{proposition}
Set
\begin{equation}
A'(X,Y,Z)=\tr_g(R(\cdot, X)A))(\cdot,Y,Z).
\end{equation}
We shall now  estimate  the function $g(A',A)$ from below. We have
\begin{equation}
\begin{array}{rcl}
&&g(A',A)= \sum_{i,j,k,l}\langle(R(e_i,e_k)A)(e_i,
e_j,e_l),A(e_k,e_j,e_l)\rangle\\
&&= -\sum_{i,j,k,l}[A(R(e_i,e_k)e_i, e_j,e_l)A(e_k,e_j,e_l)\\
&&\ \ \ \ \ \ \ \ \ \ \ \ \ \ \ \ +A(e_i, R(e_i,e_k)e_j,
e_l )A(e_k,e_j,e_l)]\\
&&\ \ \ \ \  \ \ \ \ \ \ \ \ \ \ \ \ -\sum_{i,j,k,l}A(e_i, e_j,
R(e_i,e_k) e_l)A(e_k,e_j,e_l).
\end{array}
\end{equation}
In the last term we  interchange the indices $j$ and $l$. Since $A$
is symmetric, we get
\begin{equation}
\begin{array}{rcl}
&&g(A',A)= -\sum_{i,j,k,l}[A(R(e_i,e_k)e_i,
e_j,e_l)A(e_k,e_j,e_l)\\
&&\ \ \ \ \ \ \ \ \ \ \ \ \ \ \ \ \ \ \ \ \ \ \ \ +2A(e_i,
R(e_i,e_k)e_j, e_l )A(e_k,e_j,e_l)].
\end{array}
\end{equation}
For a fixed index $l$ we have
\begin{eqnarray*}
&&-\sum_{i,j,k}[A(R(e_i,e_k)e_i,e_j,e_l)A(e_k,e_j,e_l)
+2A(e_i,R(e_i,e_k)e_j, e_l )A(e_k,e_j,e_l)]\\
&&=-\sum_{i,j,k}[A_{e_l}(R(e_i,e_k)e_i,e_j)A_{e_l}(e_k,e_j)
 +2A_{e_l}(e_i, R(e_i,e_k)e_j )A_{e_l}(e_k,e_j)]
\end{eqnarray*}
By Lemma \ref{crucial}  we now obtain
\begin{equation}\label{g(A',A)}
g(A',A)\ge (n+1)\psi H_3.
\end{equation}
By (\ref{Simons}), (\ref{g(A',A)}), (\ref{deltaA} and Proposition
\ref{nomizu} we get
\begin{equation}\label{deltapsi}
\Delta \psi\ge 2(n+1)\psi H_3+\frac{2(n+1)}{n(n-1)}\psi^2
\end{equation}

We shall now cite a theorem on weak solutions of differential
inequalities for the Laplacian of non-negative functions. The
following version of this theorem, proved  in \cite{LSZ}, is
sufficient for our purposes

\begin{thm}\label{yau}
Let $(M,g)$ be a complete Riemannian manifold with Ricci tensor
bounded from below. Suppose that $\psi$ is a non-negative continuous
function and a weak solution of the differential inequality
\begin{equation}
\Delta \psi \ge b_0\psi ^k-b_1\psi^{k-1}-...-b_{k-1}\psi -b_k,
\end{equation}
where $k>1$ is an integer and $b_0>0$, $b_1\ge 0$,..., $b_k\ge 0$.
Let $N$ be the largest root of the polynomial equation
\begin{equation}
b_0\psi ^k-b_1\psi^{k-1}-...-b_{k-1}\psi -b_k=0.
\end{equation}
Then
\begin{equation}
\psi(p)\le N
\end{equation}
for all $p\in M$.
\end{thm}
 We have, see (\ref{for R=oR}),
 \begin{equation}\label{frombelow}
 \widehat\Ric\ge \Ric\ge
(n-1)H_2,\end{equation}
 that is, $\widehat\Ric$ is bounded from
below. Since $H_3\le 0$, by Theorem \ref{yau} and (\ref{deltapsi})
we have
\begin{equation}\label{main}
\psi\le -n(n-1)H_3.
\end{equation}
Let $X$ be a unit vector. Using  (\ref{Ric+oRic}) we  now obtain
\begin{eqnarray*}
&& \widehat {\Ric}(X,X)=\Ric (X,X)+g(K_X,K_X)
\le\Ric (X,X) + g(K,K)=\Ric(X,X)+\psi\\
&&\ \ \ \ \ \ \ \ \ \ \ \ \ \ \le (n-1)H_1-n(n-1)H_3=-(n-1)^2H_3
+\frac{(n-1)n}{2}\varepsilon .
\end{eqnarray*}
Combining this with (\ref{frombelow}) one gets the following
estimation of the Ricci tensor $\widehat \Ric$
\begin{equation}
(n-1)H_3+\frac{(n-1)(n-2)}{2}\varepsilon\le\widehat\Ric\le-(n-1)^2H_3
+\frac{(n-1)n}{2}\varepsilon.
\end{equation}
In order to estimate the scalar curvature $\hat\rho$ we use
(\ref{theorema_egregium}) and (\ref{main}). We get
\begin{equation}
n(n-1)H_3+\frac{(n-1)(n-2)n}{2}\varepsilon\le\hat\rho \le
n(n-1)(H_1-H_3)=\frac{n^2(n-1)}{2}\varepsilon.
\end{equation}
The proof of Theorem \ref{epsilonfunkcja} is completed.\koniec

Theorem \ref{epsilonfunkcja} can be obviously formulated as follows

\begin{thm}\label{epsilonfunkcja1}
Let $(g,\nabla)$ be a trace-free conjugate symmetric statistical
structure on an $n$-dimensional  manifold $M$. Assume that $(M,g)$
is complete and the sectional $\nabla$-curvature   $k$ satisfies the
inequality {\rm(\ref{warunek12})} on $M$,  where $H_1=
H_3+\frac{n}{2}\varepsilon$, $H_2=H_1-\varepsilon$, $H_3$ is a
non-positive number and $\varepsilon$ is a non-negative function on
$M$. Then the Ricci tensor $\widehat \Ric$ of $g$  satisfies the
inequalities
\begin{equation}
(n-1)H_2\le\widehat\Ric\le(n-1)\left[(1-n)H_1+\frac{n^2}{2}\varepsilon\right].
\end{equation}
The scalar curvature $\hat\rho$ of $g$ satisfies the inequalities
\begin{equation}
n(n-1)H_2\le\hat\rho\le\frac{n^2(n-1)}{2}\varepsilon.
\end{equation}
\end{thm}
\medskip

\begin{remark}{\rm
The estimation of the Ricci tensor $\widehat\Ric$ from below in the
above theorems is  easy and it follows from (\ref{Riccihat_ge}). The
estimation of the Ricci tensor $\widehat\Ric$
 from above
 is not optimal in Theorems \ref{epsilonfunkcja}, \ref{epsilonfunkcja1}.
  Namely, in the case of a hyperbolic sphere, that is, in the case where
  $H_1=H_2=H_3<0$, Theorem  \ref{epsilonfunkcja} gives the estimation $\widehat \Ric\le
  -(n-1)^2H_3$ (it should be $\widehat \Ric\le 0$).
  The estimation of the scalar curvature in Theorems \ref{inIntroduction},
\ref{epsilonfunkcja1} is optimal and, in the above proof, it is not
deduced from the estimation of the Ricci tensor. }\end{remark}

\bigskip
\section{Conjugate-symmetric trace-free statistical structures with non-negative
 sectional $\nabla$-curvature}

We shall prove
\begin{thm}\label{thmpositive}
Let $(M,g)$ be a complete Riemannian manifold with  a conjugate
symmetric trace-free statistical structure  $(g,\nabla)$. If the
sectional $\nabla$-curvature is non-negative on $M$ then the
statistical structure is trivial, i.e. $\nabla=\hat\nabla$.
\end{thm}
This theorem can be deduced from the considerations of the previous
section,
 but it can be proved in an easier way, as it is shown below.
Namely, consider  the non-negative function $\varphi$ on $M$ given
by
\begin{equation}
\varphi_x={\rm max}_{U\in \mathcal U_x}A(U,U,U),
\end{equation}
where $\mathcal U_x$ is the unit hypersphere in $T_xM$, $x\in M$.
The function $\varphi$ is  continuous on $M$. Let $p\in M$ be a
fixed point and $V\in \mathcal U_p$ be a vector for which $A(U,U,U)$
attains its maximum on $\mathcal U_p$. One observes (see, e.g.
\cite{BW4} the proof of Theorem 5.6) that $V$ is an eigenvector of
$K_V$ and if $e_1=V, e_2,..., e_n$ is an orthonormal eigenbasis of
$K_V$ with corresponding eigenvalues $\lambda_1,...,\lambda_n$ then
\begin{equation}
\lambda_1-2\lambda_i\ge 0
\end{equation}
for $i=2,...,n$. Extend $V=e_1$ and $e_2,..., e_n$  by
$\hat\nabla$-parallel transport along $\hat\nabla$-geodesics
starting at $p$. We obtain a smooth orthonormal frame field. Denote
the vector  fields again by $V=e_1$ $e_2,...,e_n$. Then we have at
$p$
\begin{equation}\label{V}
\hat \nabla {e_i}=0,\ \ \ \ \  \hat\nabla _{e_i}\hat\nabla _{e_i}V=0
\end{equation}
for $i,j=1,...,n$. Denote by $\Phi$ the function $A(V,V,V)$  defined
in a neighborhood of $p$. Of course, $\Phi_p=\varphi_p$. We have at
$p$
\begin{equation}\label{laplasjanPhi}
\Delta
\Phi=\sum_{i=1}^n(\hat\nabla_{e_i}(\hat\nabla_{e_i}A))(V,V,V).
\end{equation}
Indeed, we have
\begin{eqnarray*}
(\hat\nabla d\Phi )(X,Y)&&=X(d\Phi(Y))-d\Phi(\hat\nabla_XY)\\
&&=X[(\hat\nabla _YA)(V,V,V)+3A(\hnabla_YV,V,V)]-d\Phi(\hnabla_XY)\\
&&=
(\hnabla_X(\hnabla_YA))(V,V,V)+3(\hnabla_YA)(\hnabla_XV,V,V)+3(\hnabla_XA)(\hnabla_YV,V,V)\\
&&\ \ \ \ \ +3A(\hnabla_X\hnabla_YV,V,V)
+6A(\hnabla_YV,\hnabla_XV,V)-d\Phi (\hnabla_XY).
\end{eqnarray*}
Thus, by (\ref{V}), we get (\ref{laplasjanPhi}) at $p$. We now have
at $p$
\begin{eqnarray*}
&&\Delta\Phi=\sum_{i=1}^n\hnabla_{e_i}((\hnabla_{e_i}A)(V,V,V))
=\sum_{i=1}^n\hnabla_{e_i}((\hnabla_VA)(e_i,V,V)\\
&& =\sum_{i=1}^n(\hnabla_{e_i}(\hnabla_VA))(e_i,V,V)\\
&&=\sum_{i=1}^n(\hat R(e_i,V)A)(e_i,V,V)+\sum_{i=1}^n(\hnabla_V(\hnabla_{e_i}A))(e_i,V,V)\\
&&=\sum_{i=1}^n(\hat
R(e_i,V)A)(e_i,V,V)+\sum_{i=1}^n\hnabla_V((\hnabla_{e_i}A)(e_i,V,V))\\
&&=\sum_{i=1}^n(\hat R(e_i,V)A)(e_i,V,V).
\end{eqnarray*}
In the last computations we  used both assumptions: the conjugate
symmetry and  the trace-freeness of the statistical structure. By a
straightforward computation one also gets at $p$
\begin{equation}
-\sum_{i=1}^n([K_{e_i},
K_V]A)(e_i,V,V)=\sum_{i=1}^n\lambda_i^2(3\lambda_1-2\lambda_i)
\end{equation}
and
\begin{equation}
\sum_{i=1}^n(R(e_i,V)A)(e_i,V,V)=\sum_{i=1}^n(\lambda_1-2\lambda_i)k_{i1}.
\end{equation}
Assume now that the sectional $\nabla$-curvature is bounded from
below by a number $N$. Using the equality  $\hat R=R-[K,K]$ and the
relations $\lambda_1-2\lambda_i\ge 0$, $\Phi=\lambda_1\ge 0$,
$\Phi=\varphi$ at $p$, we get at $p$
\begin{equation}
\begin{array}{rcl}
&&\Delta\Phi=\sum_{i=1}^n (\lambda_1-2\lambda_i)k_{1i}+\lambda_1^3+
\sum_{i=2}^n\lambda_i^2(3\lambda_1-2\lambda_i)\\
&&\ge \sum_{i=2}^n(\lambda_1-2\lambda_i)N+\Phi ^3=
(n+1)N\Phi+\Phi^3.
\end{array}
\end{equation}
It follows that  the function $\varphi$ is a weak solution of the
differential  inequality
\begin{equation}
\Delta\varphi\ge (n+1)N\varphi+\varphi^3.
\end{equation}
Since $\widehat{\Ric} $ is clearly bounded from below,
 by Theorem \ref{yau} we obtain that if $N\le 0$ then

\begin{equation}
\varphi(x)\le \sqrt{-(n+1)N}
\end{equation}
for all $x\in M$. If $N=0$ we get $\varphi\equiv 0$ which means that
$K\equiv 0$. Theorem \ref{thmpositive} is proved. \koniec

We also proved
\begin{proposition}\label{do_complete}
Let $(M,g)$ be a complete Riemannian manifold and $(g,\nabla)$ a
trace-free conjugate symmetric statistical structure on $M$. If the
sectional $\nabla$-curvature is bounded from below by a non-positive
number $N$ then for any unit tangent vector $U\in TM$ we have
\begin{equation}
A(U,U,U)\le \sqrt{-(n+1)N}.
\end{equation}
\end{proposition}

\section{Proof of Theorem \ref{inIntroduction}}

We shall now prove  Theorem \ref{inIntroduction}. Assume that the
statistical  sectional curvature is bounded from below and above,
that is,   the inequalities
\begin{equation}
H_2\le k(\pi)\le H_1
\end{equation}
are satisfied, where $H_1, H_2$ are real numbers. If $H_2< 0$ then
$H_3=H_2-\frac{n-2}{2}\varepsilon< 0$ and we  can use Theorem
\ref{epsilonfunkcja} to get the first assertion of Theorem
\ref{inIntroduction}. If $H_2\ge 0$ then we can use Theorem
\ref{thmpositive}. The fact that the Ricci tensor of $g$ is bounded
trivially follows from the fact that the ordinary sectional
curvature of $g$ is equal to the sectional $\nabla$-curvature. If
$H_2>0$ then  $\widehat\Ric\ge (n-1)H_2>0$. By Myers' theorem, $M$
is compact and its first fundamental group is finite. This completes
the proof of Theorem \ref{inIntroduction}.


\bigskip

\bigskip

\begin{thebibliography}{20}


\bibitem{F} Furuhata H.,  Hasegawa I., \emph{ Submanifold theory in holomorphic statistical
manifolds} in \emph{ Geometry of Cauchy-Riemann Submanifolds; S.
Dragomir, M.H. Shahid and F.R. Al-Solamy (eds.)}, Springer
Singapore, ISBN; 978-981-10-0915-0, (2016)



\bibitem{L}  Lauritzen  S.L., \emph{Statistical manifolds},  IMS
Lecture Notes-Monograph Series 10 (1987) 163-216,

\bibitem{LSZ} Li A.-M., Simon U., Zhao G.,  \emph{Global Affine Differential Geometry of Hypersurfaces},
W. de Greuter, Berlin-New York, 1993


\bibitem{No} Noguchi M., \emph{Geometry of statistical structures},
Diff. Geom. Appl., 2 (1992), 197-222

\bibitem{NS} Nomizu K., Sasaki T., \emph{Affine Differential
Geometry, Geometry of Affine Immersions}, Cambridge University
Press, 1994



\bibitem{BW4} Opozda B., \emph{Bochner's technique for statistical
structures}, Ann. Glob. Anal. Geom., 48 (2015) 357-395


\end{thebibliography}
\end{document}